# Correcting for selection bias via cross-validation in the classification of microarray data[*]

## G. J. McLachlan[1], J. Chevelu[2] and J. Zhu[1]

*University of Queensland and University of Rennes*

**Abstract:** There is increasing interest in the use of diagnostic rules based on microarray data. These rules are formed by considering the expression levels of thousands of genes in tissue samples taken on patients of known classification with respect to a number of classes, representing, say, disease status or treatment strategy. As the final versions of these rules are usually based on a small subset of the available genes, there is a selection bias that has to be corrected for in the estimation of the associated error rates. We consider the problem using cross-validation. In particular, we present explicit formulae that are useful in explaining the layers of validation that have to be performed in order to avoid improperly cross-validated estimates.

## 1. Introduction

Microarray experiments were first described in the mid-1990s as a means of measuring the expression levels of thousands of genes simultaneously (Lipshutz et al. [4] and Schena et al. [8]). They were quickly adopted by the research community for the study of a wide range of problems in the biological and medical sciences (Quackenbush [7]). One of the most important emerging clinical applications of microarray technology concerns the diagnosis of diseases, particularly in the context of cancer diagnosis. For example, accurate diagnosis of breast cancer can spare a significant number of breast cancer patients from receiving unnecessary adjuvant systemic treatment.

There are two basic approaches to generating microarray data. In a two-colour array, two samples of ribonucleic acid (RNA), each labelled with a different dye are simultaneously hybridized to the array. Cyanine 3 (Cy3) and cyanine 5 (Cy5) are fluorescent dyes that are commonly used. The sample of interest, for example, a breast cancer sample, is labelled with one dye (say Cy5), and a reference sample, for example, normal breast tissue, is labelled with another dye (say, Cy3). With this approach, the level of gene expression is estimated by the logarithm of RNA in the sample of interest to that in the (control) reference sample. For single-colour arrays, such as the Affymetrix GeneChip, each sample is labelled and individually

---

[*]Supported by the Australian Research Council.
[1]Department of Mathematics and Institute for Molecular Bioscience, University of Queensland, Australia, e-mail: gjm@maths.uq.edu.au; justin.zhu@dsto.defence.gov.au
[2]IFSIC, Université Rennes 1 and Ecole Normale Superior Britany, France, e-mail: jonathan.chevelu@eleves.bretagne.ens-cachan.fr

incubated with an array. The level of expression for each gene is represented by a single fluorescence intensity.

There is increasing interest in the use of diagnostic rules based on microarray data. In this context, there are available tissue samples of known classification with respect to a number of classes representing various conditions; for example, presence or absence of a disease or different types of treatments. The aim is to form a prediction rule based on the gene-signature vectors of these classified tissue samples. As these samples contain the expression levels on thousands of genes, some form of variable (gene) selection is usually employed before or during the formation of the prediction rule. As the final form of the prediction rule adopted will typically depend on a much smaller subset of the genes relative to the full set, care must be exercised in estimating the error rates to ensure that the consequent selection bias is corrected for. This selection bias is often overlooked in the estimation of the error rates in the bioinformatics literature so that it is common to read of a rule based of only a few genes with a zero or near-zero assessed error rate.

We focus on a nonparametric prediction rule, namely the support machine machine and the nonparametric approach of cross-validation for the estimation of its error rates. In particular, we present explicit formulae that show and explain the extra layers of validation that need to be carried out to correct satisfactorily for the selection bias or biases in those situations where there is more than one source.

## 2. Notation

### *2.1. Prediction rules*

Before we proceed to discuss the selection bias problem in estimating the accuracy of a prediction rule, we need to introduce some notation. Also, we need to define formally the prediction rule to be considered and the method of error-rate estimation to be adopted.

Although biological experiments vary considerably in their design, the data generated by microarray experiments can be viewed as a matrix of expression levels. For $M$ microarray experiments (corresponding to $M$ tissue samples), where we measure the expression levels of $N$ genes in each experiment, the results can be represented by an $N \times M$ matrix. For each tissue, we can consider the expression levels of the $N$ genes, called its *expression signature*. Conversely, for each gene, we can consider its expression levels across the different tissue samples, called its *expression profile*. The $M$ tissue samples might correspond to each of $M$ different patients.

In the present context, the problem is to construct a prediction (discriminant) rule $r(\boldsymbol{y})$ that can accurately predict the class of origin of a tumour tissue with feature vector (gene-signature vector) $\boldsymbol{y}$, which is unclassified with respect to a known number $g \, (\geq 2)$ of distinct tissue classes, denoted here by $C_1, \ldots, C_g$. Here the signature vector $\boldsymbol{y}$ contains the expression levels of a very large number $p$ of genes. If $r(\boldsymbol{y}) = i$, then it implies that $\boldsymbol{y}$ should be assigned to the $i$th class $C_i$ ($i = 1, \ldots, g$). In applications concerned with the diagnosis of cancer, one class ($C_1$) may correspond to cancer and the other ($C_2$) to benign tumours. In applications concerned with patient survival following treatment for cancer, one class ($C_1$) may correspond to the good-prognosis class and the other ($C_2$) to the poor-prognosis class. Also, there is interest in the identification of "marker" genes that characterize the different tissue classes. This is the feature selection problem. In the situation where the intention is limited to making an outright assignment to one of the possible classes, it is perhaps more appropriate to use the term prediction rather than discriminant to describe the rule. However, we shall use either nomenclature regardless of the underlying situation. In the pattern recognition jargon, such a rule is referred to as a classifier.

In order to train the prediction rule, there are available training data $t$ consisting of $n$ tissue samples of known classification. These data are obtained from $n$ microarrays, where the $j$th microarray experiment gives the expression levels of the $p$ genes in the $j$th tissue sample $y_j$ of the training set. The vector

$$(2.1) \qquad t = (y_1^T, z_1^T, \ldots, y_n^T, z_n^T)^T,$$

denotes the training data, where

$$z_j = (z_{1j}, \ldots, z_{gj})^T$$

is the class-indicator vector, and $z_{ij}$ is one or zero according as $y_j$ comes from the $i$th class $C_i$ or not ($i = 1, \ldots, g$; $j = 1, \ldots, n$). We shall write the sample rule formed from the training data $t$ as $r(y; t)$ to show its dependence on the training data $t$.

## 2.2. Different types of error rates

For a given realization $t$ of the training data $T$, it is the conditional or actual allocation rates of a sample prediction rule $r(y; t)$ that are of central interest. They are given by

$$(2.2) \qquad ec_{ij} = \mathrm{pr}\{r(Y; t) = j \mid Y \in C_i, t\} \qquad (i, j = 1, \ldots, g).$$

That is, $ec_{ij}$ is the probability, conditional on $t$, that a randomly chosen observation from $C_i$ is assigned to $C_j$ by $r(y; t)$.

The unconditional or expected allocation rates of $r(y; t)$ are given by

$$\begin{aligned} eu_{ij} &= \mathrm{pr}\{r(Y; T) = j \mid Y \in C_i\} \\ &= E\{ec_{ij}\} \qquad (i, j = 1, \ldots, g). \end{aligned}$$

The unconditional rates are useful in providing a guide to the performance of the rule before it is actually formed from the training data.

Concerning the error rates specific to a class, the conditional probability of misallocating a randomly chosen member from $C_i$ is

$$ec_i = \sum_{j \neq i}^{g} ec_{ij} \qquad (i = 1, \ldots, g).$$

The overall conditional error rate for an entity drawn randomly from a mixture $G$ of $C_1, \ldots, C_g$ in proportions $\pi_1, \ldots, \pi_g$, respectively, is

$$ec = \sum_{i=1}^{g} \pi_i ec_i.$$

The individual class and overall unconditional error rates, $eu_i$ and $eu$, are defined similarly. It is noted that the overall error rate may give a misleading summary of the error rates when the class-sample sizes are disparate. However, it is often reasonable to take the prior probabilities to be the same, where the latter are now interpreted as the class-prior probabilities adjusted (multiplicatively) by the relative importance of the costs of misallocation. For example, in the case of two classes, the cost of misallocation is often much greater for the rarer class; see McLachlan [5], Page 9.

If $r(\boldsymbol{y};\boldsymbol{t})$ is constructed from $\boldsymbol{t}$ in a consistent manner with respect to the Bayes rule $r_o(\boldsymbol{y})$, then

$$\lim_{n\to\infty} eu = e_o,$$

where $e_o$ denotes the optimal error rate. Interest in the optimal error rate in practice is limited to the extent that it represents the error of the best obtainable version of the sample-based rule $r(\boldsymbol{y};\boldsymbol{t})$. The Bayes or optimal rule $r_o(\boldsymbol{y})$ is defined by

$$(2.3) \qquad r_o(\boldsymbol{y}) = \arg\max_i \tau_i(\boldsymbol{y}),$$

where

$$\begin{aligned} \tau_i(\boldsymbol{y}) &= \mathrm{pr}\{\boldsymbol{Y}\in C_i \mid \boldsymbol{y}\} \\ &= \pi_i f_i(\boldsymbol{y})/f(\boldsymbol{y}) \end{aligned} \qquad (2.4)$$

is the posterior probability that $\boldsymbol{y}$ belongs to the $i$th class $C_i$ ($i=1,\ldots,g$). Here $f_i(\boldsymbol{y})$ denotes the density of $\boldsymbol{y}$ in class $C_i$ and

$$(2.5) \qquad f(\boldsymbol{y}) = \sum_{i=1}^{g} \pi_i f_i(\boldsymbol{y})$$

is the unconditional density of $\boldsymbol{y}$.

## 3. Support vector machine

### 3.1. Definition

In the sequel, we shall focus on the use of a nonparametric classifier, namely the support vector machine (SVM), as introduced by Vapnik [12]. Advantages of an SVM in the present context, where the number of feature variables (genes) $p$ is so large relative to the sample size $n$, are that it is able to be fitted to all the genes and that its performance appears not to be too affected by using the full set of genes. However, in practice, some form of gene selection would generally be contemplated. Another advantage of the SVM (with a linear kernel) is that gene selection can be undertaken fairly simply using the vector of weights as the criterion.

For an SVM with linear kernel, the rule $r(\boldsymbol{y};\boldsymbol{t})$ in the case of $g=2$ classes is equal to 1 or 2, according as the sign of

$$(3.1) \qquad \hat{\beta}_0 + \hat{\boldsymbol{\beta}}^T \boldsymbol{y}$$

is positive or negative. In (3.1), $\hat{\beta}_v = (\hat{\boldsymbol{\beta}})_v$ denotes the (fitted) coefficient of the expression level $y_v$ for gene $v$.

The SVM learning algorithm (with linear kernel) aims to find the separating hyperplane (3.1) that is maximally distant from the training data of the two classes. When the classes are linearly separable, the hyperplane is located so that it has maximal margin (that is, so that there is maximal distance between the hyperplane and the nearest point in any of the classes) When the data are not separable, there is no separating hyperplane; in this case, the aim is still to try to maximize the margin but allow some classification errors subject to the constraint that the total error (distance from the hyperplane on the wrong side) is less than a constant (Vapnik [12]).

### 3.2. Recursive feature elimination (RFE)

As the number of genes is much greater than the number of tissues, consideration is usually given initially to reducing the number of genes. Frequently, this is undertaken in some ad hoc manner before a more formal method of feature selection is adopted in conjunction with the choice of prediction rule. For example, one such commonly used method in the case of two classes is to rank the features on the basis of the magnitude of the (pooled) two-sample $t$-test.

For a SVM with a linear kernel, Guyon et al. [3] have shown that a good guide to the relative importance of the variables (genes) is given by the relative size of the absolute values of their fitted coefficients $\hat{\beta}_v$ (that is, the weights). Hence a ranking of the discriminatory power of the genes can be given by ranking the genes from top to bottom on the basis of the absolute values of the weights $\hat{\beta}_v$. This is what called a wrapping method, as the gene reduction method is embedded in the prediction rule.

In applying the SVM in this study, we adopt the selection procedure of Guyon et al. [3], who used a backward selection procedure, which they termed recursive feature elimination (RFE). It considers initially all the available genes, which are ranked according to their weights and the bottom-ranked genes discarded. The SVM is then refitted to the remaining genes, which are then reranked according to their new weights. Again, the bottom-ranked genes are discarded, and so on. In the applications to follow on microarray data, we first discarded enough bottom-ranked genes so that the number retained was the greatest power of 2 (less than the original number of genes). We then proceeded sequentially to discard half the current number of genes on each subsequent step. Initially, the error rate usually falls as genes are deleted, but generally, it will start to rise once a sufficiently large number of genes have been deleted. As noted by Guyon et al. [3], the process can be refined at the expense of greater computational time by deleting fewer variables on each step. One suggested procedure is to apply RFE by removing chunks of features in the first few iterations and then removing one feature at a time once the feature set size reaches a few hundred.

### 4. Estimates of error rates based on cross-validation

In practice, the error rates of a prediction rule $r(\boldsymbol{y}; \boldsymbol{t})$ are unknown as they depend on the unknown class-conditional distributions. Hence they must be estimated. In the above, we have used the notation $r(\boldsymbol{y}; \boldsymbol{t})$ to denote a prediction rule based on the training data $\boldsymbol{t}$. It is implicitly assumed that all the available genes $p$ are being used in the formation and application of this rule. In the case where only a subset of the genes are used, we write the rule as $r\{y; \boldsymbol{t}, s_d(\boldsymbol{t})\}$ to denote that it is formed using the subset $s_d(\boldsymbol{t})$. The latter denotes a subset of $d$ genes that has been selected by the adopted method of gene selection employed on the training data $\boldsymbol{t}$. Here with the SVM rule, we are using recursive feature elimination (RFE).

The apparent error rate of $r\{\boldsymbol{y}; \boldsymbol{t}, s_d(\boldsymbol{t})\}$ is given by the proportion of the tissue samples misallocated when this rule is applied to the training data $\boldsymbol{t}$. It can be expressed therefore as

$$(4.1) \qquad A\{s_d(\boldsymbol{t})\} = \frac{1}{n} \sum_{i=1}^{g} \sum_{j=1}^{n} z_{ij} Q[i, r\{\boldsymbol{y}_j; \boldsymbol{t}, s_d(\boldsymbol{t})\}]$$

where $Q[u, v]$ is one if $u \neq v$ and zero if $u = v$, and $z_{ij} = (z_j)_i$, $i = 1, \ldots, g; j = 1, \ldots, n$. Unless the sample size $n$ is large relative to the number of genes $d$ being used, it will provide too optimistic an assessment of the error rate. In the sequel, we focus on the use of cross-validation to correct the apparent error rate for its downward bias.

The leave-one-out or $n$-fold cross-validated rate is given by

$$(4.2) \qquad A^{n\mathrm{CV}}\{s_d(\boldsymbol{t})\} = \frac{1}{n} \sum_{i=1}^{g} \sum_{j=1}^{n} z_{ij} Q[i, r\{\boldsymbol{y}_j; \boldsymbol{t}_{(j)}, s_d(\boldsymbol{t}_{(j)})\}]$$

where $\boldsymbol{t}_{(j)}$ denotes the training data with the $j$th expression signature vector $\boldsymbol{y}_j$ and its class label $\boldsymbol{z}_j$ deleted; that is, with $(\boldsymbol{y}_j^T, \boldsymbol{z}_j^T)^T$ deleted. As the notation implies, we have to select a new subset of genes, $s_d(\boldsymbol{t}_{(j)})$, for the training subset $\boldsymbol{t}_{(j)}$ used on the $j$th validation trial ($j = 1, \ldots, n$). We shall refer to the cross-validated error rate (4.2) as the external cross-validated rate as on each validation trial, gene selection is undertaken externally on each training subset $\boldsymbol{t}_{(j)}$.

If we were to use the same subset $s_d(\boldsymbol{t})$ of $d$ genes as obtained originally from the full training data set $\boldsymbol{t}$ on each validation trial, then there would be a selection bias as illustrated, for example, in Ambroise and McLachlan [1]. We can express this ordinary or internal cross-validated error rate as

$$(4.3) \qquad A^{n\mathrm{CVI}}\{s_d(\boldsymbol{t})\} = \frac{1}{n} \sum_{i=1}^{g} \sum_{j=1}^{n} z_{ij} Q[i, r\{\boldsymbol{y}_j; \boldsymbol{t}_{(j)}, s_d(\boldsymbol{t})\}].$$

We use the term "ordinary" or "internal" here to mean that cross-validation is being used as it would be in the situation where the rule was based on $d$ genes with no selection process employed to choose the $d$ genes.

As the leave-one-out or $n$-fold CV rate is highly variable, it is common to use five- or ten-fold cross-validation. For example, with the latter, we can divide the $n$ training observations $(\boldsymbol{y}_j^T, \boldsymbol{z}_j^T)^T$ in $\boldsymbol{t}$ into 10 blocks (subsamples) $\boldsymbol{B}_1, \ldots, \boldsymbol{B}_{10}$ of approximatively equal size. We let $n_k$ denote the size of $\boldsymbol{B}_k$ ($k = 1, \ldots, 10$). In this case we can express the ten-fold cross-validated rate as

$$(4.4) \qquad A^{10\mathrm{CV}}\{s_d(\boldsymbol{t})\} = \frac{1}{n} \sum_{i=1}^{g} \sum_{k=1}^{10} \sum_{j \in \boldsymbol{B}_k} z_{ij} Q[i, r\{\boldsymbol{y}_j; \boldsymbol{t}_{\boldsymbol{B}_k}, s_d(\boldsymbol{t}_{(\boldsymbol{B}_k)})\}]$$

where $s_d(\boldsymbol{t}_{(\boldsymbol{B}_k)})$ denotes the selected subset of size $d$ based on the training data with the $k$th block $\boldsymbol{B}_k$ deleted. In equation (4.4), the third summation is over all labels $j$ of the tissue samples that belong to $\boldsymbol{B}_k$. Often this cross-validation is carried out in stratified form in that the classes are represented in each fold in approximately the same proportions as in the original training set.

Unlike $n$-fold cross-validation, the value of the estimate (4.4) with ten-fold cross-validation is not unique, as the training data does not have a unique split into $k = 10$ validation blocks. We can attempt to reduce the variance of the estimated error rate by calculating (4.4) for a number of splits of the training data $t$ into ten blocks and taking the average. This repeated cross-validation was not used in the results reported here.

## 5. Illustration of selection bias (subset of fixed size)

We demonstrate the selection bias that can occur when we estimate the error rate of a prediction rule based on a subset of the available genes via cross-validation without taking into account that the subset has been selected in some optimal way. As noted by Ambroise and McLachlan [1], this selection bias is often overlooked in the bioinformatics literature. We consider the breast cancer data set of van 't Veer et al. [11]. The data as analyzed here consist of the expression levels of 5,422 genes in tumours from 78 lymph-node negative patients with sporadic breast primary breast cancer categorized into two classes of patients $C_1$ and $C_2$, with $n_1 = 44$ in $C_1$ representing a good-prognosis class (that is, those who remained metastasis free after a period of more than five years), and with $n_2 = 34$ in a poor-prognosis class (those who developed distant metastases within five years). The patients were young (less than 55 years in age) and had lymph-node negative tumours. Further, they had not received adjuvant therapy, which is likely to modify outcome, and were diagnosed with breast cancer between 1983 and 1994, making a follow-up of 10 years or more possible.

We applied a support vector machine with RFE to these data, and the (ten-fold) cross-validated error rates at each stage of the selection procedure are displayed in Table 1. The first column gives the error rate of the (internal) cross-validated error rate $A^{10CVI}\{s_d(t)\}$ for which an external cross-validation has not been implemented. The second column gives the increased estimate $A^{10CV}\{s_d(t)\}$ using an external validation. That is, it uses the ten-fold version (4.4)

More specifically, consider the entries of 0.15 and 0.29 for $A^{(10CVI)}\{s_8(t)\}$ and $A^{(10CV)}\{s_8(t)\}$, respectively, for the SVM formed from the remaining 8 genes during the RFE process. With the former, the subset of 8 genes as obtained by RFE applied to the full data set is used on each of the ten validation trials. But with the latter, the selection procedure RFE is run on each of the ten validation trials, starting with all the genes, to obtain a possibly new reduced subset of 8 genes. The fact that we do not always obtain the same 8 genes on each on the ten validation trials can be used to identify potential marker genes.

## 6. Selection bias: Optimal subset of unrestricted size

In practice, having selected the "best" subset of a particular size, say $d$ variables (genes), attention turns to seeking the best subset over all sizes $d$. For example, on comparing the values of the estimated error rate $A^{10\text{CV}}\{s_d(t)\}$ in Table 1 for the values of $d$ considered, we can see that the minimum value of this estimate occurs for $d = 8$ and $d = 256$ genes. Thus in practice, consideration might be given to using the subset of 8 genes. But we know that its estimated error rate $A^{10\text{CV}}\{s_8(t)\}$ will not be an unbiased estimate of the error rate of the SVM with these 8 genes in its application to future tissue samples. We can correct for this additional selection bias due to the optimization of the $A^{10\text{CV}}\{s_d(t)\}$ over various sizes $d$ by performing a second layer of cross-validation.

TABLE 1
*Cross-validated error rates of SVM with RFE applied to 5,422 genes on $n_1 = 44$ good-prognosis patients ($C_1$) and $n_2 = 34$ poor-prognosis patients ($C_2$)*

| Number of genes | Internal CV error rate | External CV error rate |
|---|---|---|
| 1 | 0.28 | 0.40 |
| 2 | 0.19 | 0.40 |
| 4 | 0.21 | 0.42 |
| 8 | 0.15 | 0.29 |
| 16 | 0.18 | 0.38 |
| 32 | 0.12 | 0.38 |
| 64 | 0.10 | 0.33 |
| 128 | 0.12 | 0.32 |
| 256 | 0.17 | 0.29 |
| 512 | 0.15 | 0.31 |
| 1024 | 0.19 | 0.32 |
| 2048 | 0.22 | 0.35 |
| 4096 | 0.31 | 0.37 |
| 5422 | 0.37 | 0.37 |

More specifically, suppose that $r\{\boldsymbol{y}; \boldsymbol{t}, s^*(\boldsymbol{t})\}$ denotes the rule formed using the subset of genes $s^*(\boldsymbol{t})$ that minimizes $A^{10\text{CV}}\{s_d(\boldsymbol{t})\}$ over all values of $d$ that have been considered. That is,
$$s^*(\boldsymbol{t}) = s_h(\boldsymbol{t}),$$
where

(6.1)
$$h = \arg\min_d A^{10\text{CV}}\{s_d(\boldsymbol{t})\}.$$

In the case where equation (6.1) is satisfied by more than one value of $d$, we can work with the smallest value of $d$. It is clear that $A^{10\text{CV}}\{s_d(\boldsymbol{t})\}$ will underestimate the true error rate of $r\{\boldsymbol{y}; \boldsymbol{t}, s^*(\boldsymbol{t})\}$.

We can correct for this additional selection bias in forming our final rule by the optimization of a sequence of rules by performing an additional layer of cross-validation. More specifically, to train the rule based on $\boldsymbol{t}_{(B_k)}$ on the $k$th validation trial $(k = 1, \ldots, 10)$, we need to perform an additional cross-validation to form the optimal rule $r\{\boldsymbol{y}; \boldsymbol{t}_{(B_k)}, s^*(\boldsymbol{t}_{(B_k)})\}$, where $s^*(\boldsymbol{t}_{(B_k)})$ denotes the optimal subset over all subsets for a rule based on the training data with the $k$th block deleted. We can do this using ten-fold cross-validation; or, since $\boldsymbol{t}_{(B_k)}$ consists of 9 blocks $\boldsymbol{B}_h$ ($h = 1, \ldots, 10; h \neq k$), it is convenient to use nine-fold cross-validation. Accordingly, we can provide an approximate unbiased estimate of the error rate, given by

(6.2)
$$A^{10\text{CV}2}\{s^*(\boldsymbol{t})\} = \frac{1}{n}\sum_{i=1}^{g}\sum_{k=1}^{10}\sum_{j \in B_k} z_{ij} Q[i, r\{\boldsymbol{y}_j; \boldsymbol{t}_{(B_k)}, s^*(\boldsymbol{t}_{(B_k)})\}]$$

where
$$s^*(\boldsymbol{t}_{(B_k)}) = s_{h_k}(\boldsymbol{t}_{(B_k)})$$
and
$$h_k = \arg\min_d \sum_{i=1}^{g}\sum_{\substack{k'=1 \\ k' \neq k}}^{10}\sum_{j \in B_{k'}} \frac{z_{ij} Q[i, r\{\boldsymbol{y}_j; \boldsymbol{t}_{(B_k, B_{k'})}, s_d(\boldsymbol{t}_{(B_k, B_{k'})})\}]}{n - n_k}$$

and $\boldsymbol{t}_{(B_k, B_{k'})}$ denotes the full training data set $\boldsymbol{t}$ with the $k$th and the $k'$th blocks, $\boldsymbol{B}_k$ and $\boldsymbol{B}_{k'}$ deleted, and $n_k$ is the number of samples in the $k$th block.

Applying this estimate (6.2) to the data in Table 1, we find that the value of $A^{10\text{CV}2}\{s^*(\boldsymbol{t})\}$ is 0.37. This compares with the value of 0.29 for the estimate $A^{10\text{CV}}\{s_d(\boldsymbol{t})\}$ for $d = 8$ genes. Thus the effect of performing the second layer of cross-validation is to increase the estimate of the error rate from 0.29 to 0.37.

## 7. Additional bias from preliminary screening

In the previous sections, we have considered the SVM prediction rule for some microarray data from the breast cancer study of van 't Veer et al. [11]. It sought to distinguish between patients who had the same stage of disease but a different response to treatment and a different overall outcome. A signature vector of 70 genes was identified on the basis of the correlation between a gene and the class label, which is equivalent to using the (pooled) two-sample $t$-statistic to rank the genes. They called these 70 genes the prognostic marker genes. The superiority of this discriminant analysis approach based on this 70-gene signature as compared with traditional clinical staging is under dispute. Also, Ein-Dor et al. [2] in a study of this data set concluded that there are several sets of 70 genes with equal predictive powers, and this is supported by subsequent studies including Michiels et al. [6] and our own results which are to be reported elsewhere.

Note that there is a selection bias if a prediction rule is formed form these 70 genes without taking into account that they are the "top" 70 genes in a much larger list of genes. To illustrate this, we formed a SVM from the data set of 78 breast cancers using just these 70 genes. In Table 2, we report the value of the (internal) cross-validated error rate $A^{10\text{CVI}}\{s_d(\boldsymbol{t}^{(70)})\}$ for the values of $d$ considered, where $\boldsymbol{t}^{(70)}$ denotes the training data with expression levels available only on the "top" 70 genes. That is,

$$(7.1) \qquad A^{10\text{CVI}}\{s_d(\boldsymbol{t}^{(70)})\} = \frac{1}{n} \sum_{i=1}^{g} \sum_{k=1}^{10} \sum_{j \in \boldsymbol{B}_k} z_{ij} Q[i, r\{\boldsymbol{y}_j; \boldsymbol{t}^{(70)}_{(\boldsymbol{B}_k)}, s_d(\boldsymbol{t}^{(70)}_{(\boldsymbol{B}_k)})\}].$$

On comparing them with the results in Table 1 for the SVM using RFE starting from 5,422 genes, we can see that there is a selection bias in starting with these "top" 70 genes. In practice, we can correct for this bias provided the expression levels are available for the larger set of genes. This raises the question of how large a set of genes we should start with in order to avoid this type of bias. Zhu et al. [13] concluded that the initial set of genes has to be relatively small relative to the total number of genes available for this bias to be of practical significance.

In the second column of Table 2, we have listed the values of the (external) cross-validated error rate given by

$$(7.2) \qquad A^{10\text{CV}}\{s_d(\boldsymbol{t}^{(70)})\} = \frac{1}{n} \sum_{i=1}^{g} \sum_{k=1}^{10} \sum_{j \in \boldsymbol{B}_k} z_{ij} Q[i, r\{\boldsymbol{y}_j; \boldsymbol{t}^{(70_k)}_{(\boldsymbol{B}_k)}, s_d(\boldsymbol{t}^{(70_k)}_{(\boldsymbol{B}_k)})\}]$$

where $\boldsymbol{t}^{(70_k)}$ denotes the training data for the "top" 70 genes found when using the $k$th training subset $\boldsymbol{t}_{(\boldsymbol{B}_k)}$ with expression levels available on the full set of 5,422 genes on all tissues apart from those in the $k$th block $\boldsymbol{B}_k$, and $\boldsymbol{t}^{(70_k)}_{(\boldsymbol{B}_k)}$ denotes $\boldsymbol{t}^{(70_k)}$ with the $k$th block $\boldsymbol{B}_k$ deleted.

The point to be stressed here is that if only the subset of selected genes is available, then there is no way to correct for the selection bias in working with these "top" genes. Also, it is noted that the above approach of selecting the top genes individually via univariate methods is not the optimal multivariate method of selection for a prediction rule.

TABLE 2
*Cross-validated error rates of SVM with RFE applied to the top 70 genes on $n_1 = 44$ good-prognosis patients $(C_1)$ and $n_2 = 34$ poor-prognosis patients $(C_2)$ without and with bias correction for starting with a top subset of the genes*

| Number of genes | Internal CV error rate | External CV error Rate |
|---|---|---|
| 1 | 0.40 | 0.40 |
| 2 | 0.31 | 0.37 |
| 4 | 0.24 | 0.40 |
| 8 | 0.23 | 0.32 |
| 16 | 0.22 | 0.33 |
| 32 | 0.19 | 0.35 |
| 64 | 0.23 | 0.32 |
| 70 | 0.23 | 0.32 |

## 8. Another selection bias

The selection biases discussed in the previous work arise due to the fact that some or all of the data used to form the rule are involved in some way with the testing of the rule. One way of avoiding such biases is to use a holdout method. The available data are split into disjoint training and test sets. The prediction rule is formed from the training subset and then assessed on the test set. Clearly, this method is inefficient in its use of the data, particularly in the context of microarray data, where the number of tissue samples is usually limited.

A commonly occurring mistake in analyses of microarray data is to adopt a holdout method, but to use the test set in the selection of the genes. With this approach, the test subset plays a role in leading to the choice of the final form of the prediction rule. However, this is frequently overlooked in the bioinformatics literature (Somorjai et al. [9]). It often leads to claims that a prediction rule can be formed from only a few genes that has almost zero error rate.

## 9. The Netherlands breast cancer data

A nationwide clinical trial MINDACT (Microarray In Node Negative Disease may Avoid Chemo Therapy), is now underway in the Netherlands, in which the expression profiles for the 70-gene signature proposed by van 't Veer et al. [11] are being collected from all newly identified consenting patients with breast cancer and used as an adjunct to classic clinical staging. It can be seen from Table 1 that the data in van 't Veer et al. [11] consisting of the sporadic 78 breast cancers are of limited discriminatory value. Subsequently, van de Vijver et al. [10] considered a larger set of 295 breast cancer tissue samples, which consisted of 61 of the 78 patients in van 't Veer et al. [11]. Each of the new 234 patients, some of which were lymph-node positive, were followed for at least 5 years or to their censored time; some actually died during the followup period. There were 55 patients with censored outcomes (the occurrence or nonoccurrence of metastases within 5 years was not known).

For the 179 new patients with known outcomes combined with the 61 from the original study of van 't Veer et al. [11], we fitted a SVM with RFE, starting with all 24,481 genes. The results are displayed in Table 3. The smallest estimated error rate is 0.26 at $d = 128$ genes. If we use (6.2) to correct for the bias arising from the fact that this is the smallest over a number of subsets, we obtain an estimate of 0.28.

TABLE 3
*Cross-validated error rates starting with all available genes using predicted class labels from 61 breast cancer tumours from the van' t Veer et al. [10] study and using the true class labels from actual outcomes*

| Number of genes | External CV error rate using predicted class labels | External CV error rate using true class labels |
|---|---|---|
| 1 | 0.42 | 0.31 |
| 2 | 0.38 | 0.31 |
| 4 | 0.38 | 0.34 |
| 8 | 0.31 | 0.35 |
| 16 | 0.23 | 0.33 |
| 32 | 0.22 | 0.30 |
| 64 | 0.19 | 0.30 |
| 128 | 0.16 | 0.26 |
| 256 | 0.15 | 0.29 |
| 512 | 0.14 | 0.29 |
| 1024 | 0.15 | 0.28 |
| 2048 | 0.14 | 0.27 |
| 4096 | 0.14 | 0.27 |
| 8192 | 0.16 | 0.27 |
| 16384 | 0.17 | 0.29 |
| 24481 | 0.18 | 0.27 |

In their study, van de Vijver et al. [10] created two classes corresponding to good- and poor-prognosis by ignoring the actual outcomes where available and assigning the patients to the two classes on the basis of a rule formed from the gene-expression signatures for these 61 patients. More precisely, each of the 234 new tumours was assigned to class $C_1$ or $C_2$ according as the correlation between the elements of its gene-signature vector and the corresponding elements of $\overline{y}_1$ was greater or less than 0.4, where $\overline{y}_1$ is the mean of the gene signatures of those patients from the original 61 patients in the good-prognosis class $C_1$. Each of the 61 original tumours was reassigned to $C_1$ to $C_2$ according to whether the (cross-validated) correlation between its gene-signature vector and $\overline{y}_1$ was greater or less than 0.55 (a threshold that resulted in a 10 % rate of false negatives in the study of van 't Veer et al. [11]).

Thus if we were to use their predicted classification for these data, we would obtain results that are more optimistic than if we used the actual outcomes for the occurrence or nonoccurrence of metastases. To demonstrate this, we have listed in Table 3, the results of fitting a SVM with RFE to the 295 tissue samples with the predicted classification of van de Vijver et al. [10]. Some idea of the bias involved can be obtained by comparing them with the corresponding results in this table using the known outcomes for 240 tissue samples. Admittedly the training set is smaller for the latter, but it would not account for the differences obtained here. Further, van de Vijver et al. [10] used only the "top" 70 genes as identified from the 78 tissue samples in the van 't Veer et al. [11] study. As 61 of these 78 tissue samples are being used in the enlarged set in van de Vijver et al. [10], this will be another source of bias.

## 10. Discussion

In this study, we have described some of the ways in which a selection bias can arise in the formation of a prediction rule using a subset of the available genes selected in some optimal manner from a much larger set of genes measured on only a limited number of training samples. To illustrate these biases, we consider applications of a nonparametric rule, the support vector machine, formed on subsets of the available genes by using the feature selection procedure known as RFE (recursive feature elimination). We employ the nonparametric method of bias correction, cross-validation, to correct for these selection biases. In some situations, more than one layer of validation must be performed in order to ensure that the error rate is properly validated. Our results underscore the importance of cross-validating at all stages of the procedure used to train the prediction rule.

The methodology is demonstrated on the Netherlands breast cancer data arising from the studies of van 't Veer et al. [11] and van de Vijver et al. [10]. Identifying a gene signature for breast cancer prognosis has been a central goal in these and other microarray studies. These data have attracted considerable attention in the bioinformatics and medical literature, as their analyses have raised a number of statistical issues including aspects of the topic of selection bias as discussed here.

In van 't Veer et al. [11], a 70 gene signature (known as the Amsterdam signature) was derived from a cohort of 78 breast cancer patients. Our results here show that the SVM has an accuracy of only 0.63%. van de Vijver et al. [10] subsequently attempted to provide further validation of the Amsterdam signature by considering a larger set of 295 breast cancer patients. The results reported for the SVM trained on the breast cancers of known outcomes in this expanded study of van de Vijver et al. [10] show that it has an accuracy of approximately 72% in predicting prognosis status. This is better than the 63% obtained on the original set of 78 breast cancer tumours, but is still not high. Clearly, the molecular classification of breast cancer is still a work in progress.

tion and discovery using gene microarray and proteomics mass spectroscopy data: Curses, caveats, cautions. *Bioinformatics* **19** 1484–1491.